\documentclass[12pt,oneside,reqno]{amsart}
\usepackage{amsmath}
\usepackage{graphicx}
\usepackage{mathrsfs}
\usepackage{stmaryrd}
\usepackage{amsfonts}
\usepackage{enumerate,amsmath,amssymb,amsthm}
\newcommand{\be}{\begin{eqnarray}}
\newcommand{\ee}{\end{eqnarray}}
\newcommand{\ce}{\begin{eqnarray*}}
\newcommand{\de}{\end{eqnarray*}}
\newtheorem{theorem}{Theorem}[section]
\newtheorem{lemma}[theorem]{Lemma}
\newtheorem{remark}[theorem]{Remark}
\newtheorem{definition}[theorem]{Definition}
\newtheorem{proposition}[theorem]{Proposition}
\newtheorem{Examples}[theorem]{Example}
\newtheorem{corollary}[theorem]{Corollary}

\def\eps{\varepsilon}

\def\a{\alpha}

\def\p{\partial}
\def\d{\delta}

\def\[{{\Big[}}
\def\]{{\Big]}}
\def\<{{\langle}}
\def\>{{\rangle}}
\def\({{\Big(}}
\def\){{\Big)}}

\def\bx{{\mathbf{x}}}
\def\tr{{\rm tr}}

\def\dif{{\mathord{{\rm d}}}}

\def\no{\nonumber}
\def\bt{\begin{theorem}}
\def\et{\end{theorem}}
\def\bl{\begin{lemma}}
\def\el{\end{lemma}}
\def\br{\begin{remark}}
\def\er{\end{remark}}
\def\bd{\begin{definition}}
\def\ed{\end{definition}}
\def\bp{\begin{proposition}}
\def\ep{\end{proposition}}
\def\bc{\begin{corollary}}
\def\ec{\end{corollary}}
\def\bx{\begin{Examples}}
\def\ex{\end{Examples}}

\def\cF{{\mathcal F}}

\def\mE{{\mathbb E}}

\def\mG{{\mathbb G}}

\def\mI{{\mathbb I}}

\def\mN{{\mathbb N}}

\def\mP{{\mathbb P}}

\def\mR{{\mathbb R}}

\def\mT{{\mathbb T}}

\def\mX{{\mathbb X}}

\def\mZ{{\mathbb Z}}

\def\bP{{\mathbf P}}

\def\geq{\geqslant}
\def\leq{\leqslant}

\def\bW{{\mathbf W}}
\def\bP{{\mathbf P}}

\def\tr{{\mathrm t}}
\allowdisplaybreaks

\begin{document}

\title{A Stochastic Representation for Backward Incompressible Navier-Stokes Equations}

\date{}
\author{Xicheng Zhang}

\thanks{{\it Keywords: }Backward  Navier-Stokes Equation, Stochastic Representation,
Global Existence, Large Deviation}

\dedicatory{
School of Mathematics and Statistics\\
The University of New South Wales, Sydney, 2052, Australia,\\
Department of Mathematics,
Huazhong University of Science and Technology\\
Wuhan, Hubei 430074, P.R.China\\
Email: XichengZhang@gmail.com
 }

\begin{abstract}

By reversing the time variable we derive a stochastic representation for
backward incompressible Navier-Stokes equations in terms of stochastic Lagrangian paths,
which is similar to  Constantin and Iyer's forward formulations in \cite{Co-Iy}.
Using this representation, a self-contained
proof of local existence of solutions in Sobolev spaces are provided for
 incompressible Navier-Stokes equations in the whole space. In two dimensions
or large viscosity, an alternative
proof to the global existence is also given.
Moreover, a large deviation estimate for stochastic particle trajectories is presented when
the viscosity  tends to zero.
\end{abstract}

\maketitle
\rm

\section{Introduction}

The classical Navier-Stokes equations describe the evolution of velocity fields of
an incompressible fluid, and takes the following form with the external force zero:
\be
\label{NS}
\left\{
\begin{aligned}
&\p_t u+(u\cdot\nabla)u-\nu\Delta u+\nabla p=0,\ \ t\geq 0,\\
&\nabla\cdot u=0,\ \  u(0)=u_0,
\end{aligned}
\right.
\ee
where column vector field $u=(u^1,u^2,u^3)^\tr$ denotes
the velocity field, $p$ is the pressure and
$\nu$ is the kinematic viscosity. When the viscosity $\nu$ vanishes, the above
equation becomes the classical Euler equation:
\be
\label{Euler}
\left\{
\begin{aligned}
&\p_t u+(u\cdot\nabla)u+\nabla p=0,\ \ t\geq 0,\\
&\nabla\cdot u=0,\  \ u(0)=u_0,
\end{aligned}
\right.
\ee
which describes the motion of an ideal incompressible fluid.
The mathematical theory about Navier-Stokes equations and Euler equations has
been extensively studied and the existence of regularity solutions is still
a big open problem in modern PDEs.

Recently,  Constantin and  Iyer \cite{Co-Iy}
presented an elegant stochastic representation
for incompressible Navier-Stokes equations based on stochastic particle paths,
which is realized by an implicit stochastic differential equation:
the drift term is computed as the expected value of an expression involving
the stochastic flow defined by itself. More precisely,
let $(u,X)$ solve the following stochastic system:
\be
\label{Re}\left\{
\begin{aligned}
&X_t(x)=x+\int^t_0u_s(X_s(x))\dif s+\sqrt{2\nu}B_t,\ \ t\geq 0,\\
&A_t=X^{-1}_t,\\
&u_t=\mE\bP[(\nabla^\tr A_t) (u_0\circ A_t)],
\end{aligned}
\right.
\ee
where $B_t$ is a $3$-dimensional Brownian motion,  $\bP$ is the Leray-Hodge
projection onto divergence free vector fields,
and $\nabla^\tr A_t$ denotes the transpose of
Jacobi matrix $\nabla A_t$.
Then $u$ satisfies equation (\ref{NS}) with initial data $u_0$.
One of the proofs given in \cite{Co-Iy}
is based on a stochastic partial differential equation satisfied by the inverse flow $A$.
By using this representation, a self-contained proof of the existence of
local smooth solutions is provided in \cite{Iy,Iy-Th}.

Let $(u,p)$ solve (\ref{NS}). Notice that if we make the time change:
$$
\tilde u(t,x):=-u(-t,x),\ \ \tilde p(t,x)=p(-t,x)
\mbox{ for $t\leq 0$},
$$
then $\tilde u$ satisfies the following equation (called backward Navier-Stokes equation here):
\be
\label{Ba}
\left\{
\begin{aligned}
&\p_t \tilde u+(\tilde u\cdot\nabla)\tilde u+\nu\Delta \tilde u+\nabla p=0,\ \ t\leq 0, \\
&\nabla\cdot \tilde u=0,\ \  \tilde u(0)=u_0.
\end{aligned}
\right.
\ee
The purpose of the present paper is to give a slightly different representation for $\tilde u$
by using backward particle paths. More precisely, let $(\tilde u,X)$ solve
the following stochastic system
\be
\left\{
\begin{aligned}
&X_{t,s}(x)=x+\int^s_t\tilde u_r(X_{t,r}(x))\dif r+\sqrt{2\nu}(W_s-W_t),\
t\leq s\leq 0,\\
&\tilde u_t=\mE \bP[(\nabla^\tr X_{t,0}) (u_0\circ X_{t,0})] ,\ \ t\leq0.
\end{aligned}
\right.\label{Rep0}
\ee
Then $\tilde u$ satisfies backward incompressible Navier-Stokes equation (\ref{Ba}) with final value
$u_0$.
Intuitively, by reversing the time, the starting point is
changed as the end point. Hence, representation (\ref{Rep0})
is essentially the same as (\ref{Re}).
But, equation (\ref{Rep0}) is easier to be
dealt with mathematically since the unpleasant term $A=X^{-1}$ which usually incurs extra mathematical calculations
does not appear in (\ref{Rep0}).   Such a representation will be proved in
Section 2. We emphasize that representations (\ref{Re}) and (\ref{Rep0}) are
useful in numerical computations (cf. \cite{Be-Ma,Iy-Ma}). Direct calculations
shows that the second equation in (\ref{Rep0}) is equivalent to
\be
&&\tilde\omega_t=\nabla\times \tilde u_t=\mE [(\nabla^{-1} X_{t,0}) (\omega_0\circ X_{t,0})],\\
&&\tilde u_t=-\Delta^{-1}\nabla\times\tilde\omega_t,\label{Lp1}
\ee
where $\tilde\omega_t$ is the vorticity,
$\nabla^{-1} X_{t,0}$ is the inverse of Jacobi matrix $\nabla X_{t,0}$ and
(\ref{Lp1}) is exactly the Biot-Savart law.

We also mention other stochastic formulations for incompressible Navier-Stokes equations.
In \cite{Es-Ma-Pu}, a representation formula for the vorticity of three dimensional Navier-Stokes
equations was given by using  stochastic Largrangian paths, however, there is no
a self-contained proof of the existence given there. In \cite{Le-Sz}, Le Jan and Sznitman
used a backward-in-time branching process in Fourier space to
express the velocity field of a three-dimensional viscous fluid as the average
of a stochastic process, which then leads to a new existence theorem. In \cite{Bu},
basing on Girsanov's transformation and Bismut-Elworty-Li's formula, Busnello introduced
a purely probabilistic treatment to the existence of a unique global solution
for two dimensional Navier-Stokes equations, where the stretching term disappears, and
the non-linear equation obeyed by the vorticity has the form of Fokker-Planck equation.
Later on,  Busnello, Flandoli and Romito in \cite{Bu-Fl-Ro} carefully analyzed
an implicit probabilistic representation for the vorticity of
three dimensional Navier-Stokes equations, and
a local existence was given. In that paper, much attentions were
also paid on a probabilistic representation
formula for a general system of linear parabolic equations. Moreover, in \cite{Bu, Bu-Fl-Ro},
an interesting probabilistic representation for the Biot-Savart law was also given and analyzed
so that they can recover the velocity from the vorticity by probabilistic approach.
Recently,  Cipriano and Cruzeiro in \cite{Ci-Cr} described a stochastic variational principle for
two dimensional incompressible Navier-Stokes equations by using the Brownian motions on the group
of homeomorphisms on the torus. More recently, Cruzeiro and Shamarova \cite{Cr-Sh}
 established a connection between equation (\ref{Ba}) and a system of infinite dimensional forward-backward
stochastic differential equations on the group of volume-preserving diffeomorphisms of a flat torus.

This paper is organized as follows:
In Section 3, we shall give a self-contained proof of local existence in Sobolev spaces.
The proof is based on successive approximation or  fixed point method as in \cite{Iy}.
Therein,  Iyer considered
the spatially periodic case and worked in H\"older continuous function spaces.
When Sobolev spaces are considered, we have to overcome the difficulty due to
the non-closedness of Sobolev spaces under pointwise multiplications  and compositions.
Thus, it seems to be hard to exhibit the same proof in Section 3 for representation (\ref{Re})
due to the presence of $A=X^{-1}$.
A key point in the proof lies in that the flow map $x\mapsto X_{t,s}(x)$
preserves the Lebesgue measure, i.e., $\mathrm{det}(\nabla X_{t,s})=1$.

In Section 4, we shall give an alternative proof to the global existence
when the spatial dimension is two or the viscosity is large enough in any dimensions.
Such results are well known. For two dimensional Navier-Stokes and
Euler equations in the whole space, the global existence of smooth solutions are referred to
\cite{Be-Ma}. The global existence of regularity solutions for large viscosity is referred to
\cite{Ko-Ta}. Recently, Iyer \cite{Iy0} presented an alternative proof to the global existence for
small Reynolds. His proof is based on the decay of heat flows and
stochastic representation (\ref{Re}). Following \cite{Iy0}, we will give a
different proof based on Bismut's formula and representation (\ref{Rep0}).

Let $(u^\nu,X^\nu)$ denote the solution of equation (\ref{Rep0}). In Section 5,
as $\nu$ goes to zero,
an asymptotic probability estimate of $X^\nu$ in diffeomorphism group is presented by
the well known large deviation estimate for stochastic diffeomorphism flows.

\section{Stochastic Representation of Backward Incompressible Navier-Stokes Equations}

We begin with some notational conventions. Fix $d\geq 2$ and put
$$
\mN_0:=\{0\}\cup\mN,\ \ \mR_-:=(-\infty,0],\ \ \mI:=(d\times d)\mbox{-unit matrix}
$$
and
$$
l^2:=\left\{\sigma=(\sigma_k)_{k\in\mN}\in\mR^\mN: \|\sigma\|_{l^2}^2:=\sum\sigma_k^2<+\infty\right\}.
$$
For a differentiable transformation $X$ of $\mR^d$, the Jacobi matrix of $X$ is given by
\ce
\nabla X:=\left(
\begin{array}{cccc}
\p_1 X^1,&\p_2 X^1,&\cdots, &\p_d X^1\\
\p_1 X^2,&\p_2 X^2,&\cdots, &\p_d X^2\\
\cdots,&\cdots,&\cdots,&\cdots\\
\p_1 X^d,&\p_2 X^d,&\cdots,& \p_d X^d
\end{array}
\right),
\de
where $\p_i=\frac{\p}{\p x_i}$.
We use $\nabla^\tr X$ to denote the transpose of $\nabla X$.
For $k\in\mN_0$, let $C^k_b(\mR^d;\mR^d)$ denote the space of $k$-order continuous
differentiable vector fields on $\mR^d$ with the norm:
$$
\|u\|_{C^k_b}:=\sum_{|\a|\leq k}\sup_{x\in\mR^d}|D^\a u(x)|<+\infty,
$$
where $D^\a$ denotes the derivative with respect to the multi index $\a$.

Let $u\in C(\mR_-; C^3_b(\mR^d;\mR^d))$ and
$t\mapsto\sigma_t=\sigma(t)\in L^2_{loc}(\mR_-;\mR^d\times l^2)$
satisfy
\be
\sum_{k=1}^\infty\sigma^\cdot_k(t)\sigma^\cdot_k(t)=\mI.\label{Con}
\ee
Let $\{X_{t,s}(x), t\leq s\leq 0\}$ solve the following  SDE
\be
X_{t,s}(x)=x+\int^s_tu_r(X_{t,r}(x))\dif r+\sqrt{2\nu}\int^s_t\sigma_r\dif B_r,
\label{SDE}
\ee
where $B_t:=\{B_t^k, t\leq 0, k\in\mN\}$ is a sequence of independent standard  Brownian motions
on some probability space $(\Omega,\cF,\mP)$.
Thanks to (\ref{Con}), the diffusion operator associated to equation (\ref{SDE}) is given by
$$
L_tg\equiv \nu\Delta g+(u_t\cdot\nabla)g.
$$

We first prove the following result.
\bt\label{Th11}
Let $\phi:\mR^d\to\mR^d$ be a $C^2$-vector field satisfying
$$
|D^\a\phi(x)|\leq C(1+|x|^\beta),\ \ |\a|\leq 2,\ \beta>0,
$$
and $f\in C(\mR_-; C^2_b(\mR^d;\mR^d))$. Define
\ce
\Lambda_\phi(t,x)&:=&(\nabla^\tr X_{t,0} (x))\phi(X_{t,0}(x)),\\
\Lambda_f(t,x)&:=&\int^0_t  (\nabla^\tr X_{t,r}(x))f_r(X_{t,r}(x))\dif r
\de
and
$$
w_t(x):=\mE \Lambda_\phi(t,x)-\mE \Lambda_f(t,x).
$$
Then $w\in C^{1,2}((-\infty,0)\times\mR^d)$ satisfies the following backward Kolmogorov's equation:
\be
\p_t w_t+L_tw_t+(\nabla^\tr u_t)w_t=f_t, \ \ \lim_{t\uparrow 0}w_t(x)=\phi(x).\label{SA}
\ee
\et
\begin{proof}
Let $g$ be a twice continuously differentiable function satisfying
$$
|D^\a g(x)|\leq C(1+|x|^\beta),\ \ |\a|\leq 2,\ \beta>0.
$$
For $h>0$, by It\^o's formula we have
$$
\mE g(X_{t-h,t}(x))=g(x)+\mE\left[\int^t_{t-h}(L_rg)(X_{t-h,r}(x))\dif r\right].
$$
From this, it is easy to see that
\be
\frac{1}{h}\Big[\mE g(X_{t-h,t}(x))-g(x)\Big]&=&\frac{1}{h}\mE\left[\int^t_{t-h}(L_rg)
(X_{t-h,r}(x))\dif r\right]\no\\
&\to& (L_tg)(x)\ \ \mbox{as $h\to 0$}.\label{Es2}
\ee

Noticing that
$$
X_{t-h,0}(x)=X_{t,0}\circ X_{t-h,t}(x),
$$
we have
$$
\nabla X_{t-h,0}(x)=(\nabla X_{t,0})\circ X_{t-h,t}(x)\cdot\nabla X_{t-h,t}(x).
$$
Thus, by the independence of $X_{t-h,t}(x)$ with $X_{t,0}(x)$, we have
\ce
\mE \Lambda_\phi(t-h,x)&=&\mE[(\nabla^\tr X_{t-h,t}(x))\Lambda_\phi(t,X_{t-h,t}(x))]\\
 &=&\mE\Big[(\nabla^\tr X_{t-h,t}(x))((\mE \Lambda_\phi(t))\circ X_{t-h,t}(x))\Big]
\de
and
\ce
\mE \Lambda_f(t-h,x)&=&\mE\left[(\nabla^\tr X_{t-h,t}(x))\Lambda_f(t, X_{t-h,t}(x))\right]\\
&&+\mE\left[\int^t_{t-h}(\nabla^\tr X_{t-h,r}(x))f(r,X_{t-h,r}(x)) \dif r\right]\\
&=&\mE\left[(\nabla^\tr X_{t-h,t}(x))( (\mE\Lambda_f(t))\circ X_{t-h,t}(x))\right]\\
&&+\mE\left[\int^t_{t-h}(\nabla^\tr X_{t-h,r}(x))f_r(X_{t-h,r}(x))\dif r\right].
\de
Hence, we may write
\ce
w_{t-h}(x)&=&\mE\Big[(\nabla^\tr X_{t-h,t}(x)) w_t(X_{t-h,t}(x)) \Big]\\
&&-\mE\left[\int^t_{t-h}(\nabla^\tr X_{t-h,r}(x))f_r(X_{t-h,r}(x))\dif r\right]
\de
and
\ce
\frac{1}{h}(w_t(x)-w_{t-h}(x))
&=&-\frac{1}{h}\mE\Big[(\nabla^\tr X_{t-h,t}(x)-\mI)w_t(X_{t-h,t}(x))\Big]\\
&&-\frac{1}{h}\Big[\mE\left(w_t(X_{t-h,t}(x))\right)-w_t(x)\Big]\\
&&+\frac{1}{h}\mE\left[\int^t_{t-h}(\nabla^\tr X_{t-h,r}(x))f_r(X_{t-h,r}(x))\dif r\right]\\
&=:&I_1^h(t,x)+I_2^h(t,x)+I_3^h(t,x).
\de
Observing that
$$
\nabla X_{t-h,t}(x)-\mI=\int^t_{t-h}(\nabla  u_s)\circ X_{t-h,s}(x)\cdot\nabla X_{t-h,s}(x)\dif s.
$$
we deduce
$$
\lim_{h\downarrow 0}I_{1}^h(t,x)=-[(\nabla^\tr u_t) w_t](x).
$$
By  (\ref{Es2}) we have
$$
\lim_{h\downarrow 0}I_2^h(t,x)=-(L_tw_t)(x).
$$
Moreover, a simple limit procedure also gives
$$
\lim_{h\downarrow 0}I_3^h(t,x)=f_t(x).
$$
Combining the above calculations, we conclude that
$$
\lim_{h\downarrow 0}\frac{1}{h}(w_t(x)-w_{t-h}(x))=-(L_tw_t)(x)-[(\nabla^\tr u_t) w_t](x)+f_t(x).
$$
Equation (\ref{SA}) now follows (see \cite[p.124]{Fr} for more details).
\end{proof}
\br
A more general Feynman-Kac formula for a
deterministic system of parabolic equations was given in \cite{Bu-Fl-Ro}. However,
 the proof is simpler in our case. In representation
(\ref{Re}), if we define $w_t:=\mE[(\nabla^\tr A_t) (u_0\circ A_t)]$, then $w_t$ also satisfies
(\ref{SA}) with $f=0$(see \cite[p.343, (4.5)]{Co-Iy}).
\er

Basing on this theorem, we can give a stochastic representation for
backward Navier-Stokes equation (\ref{Ba}) as in \cite{Co-Iy}.
\bt
Let $\nu\geq 0$ and $u_0\in C^2_b(\mR^d;\mR^d)$ a deterministic divergence-free vector field,
and $f\in C(\mR_-; C^2_b(\mR^d;\mR^d))$.
Suppose that $\sigma$ satisfies (\ref{Con}), and  $(u,X)$ solves the stochastic system:
\be
X_{t,s}(x)&=&x+\int^s_tu_r(X_{t,r}(x))\dif r+\sqrt{2\nu}\int^s_t\sigma_r\dif B_r,\
t\leq s\leq 0,\label{Se1}\\
u_t&=&\bP\mE \Lambda^u_{u_0}(t)-\bP\mE \Lambda^u_{f}(t),\ \ t\leq0,\label{Rep}
\ee
where $\bP$ is the Leray-Hodge projection onto divergence free vector fields,
and $\Lambda^u_{u_0}$ and $\Lambda^u_{f}$ are given by
\ce
\Lambda^u_{u_0}(t,x)&:=&(\nabla^\tr X_{t,0}(x))u_0(X_{t,0}(x)), \\
\Lambda^u_{f}(t,x)&:=&\int^0_t(\nabla^\tr X_{t,r}(x))f_r(X_{t,r}(x))\dif r.
\de
Then $u$ satisfies the backward incompressible Navier-Stokes equation:
\be
\label{Back}\left\{
\begin{aligned}
&\p_t u+(u\cdot\nabla)u+\nu\Delta u+\nabla p=f,\quad  t\leq 0,\\
&\nabla\cdot u=0, \quad u(0,x)=u_0(x).
\end{aligned}
\right.
\ee
Conversely, if $u$ solves  backward Navier-Stokes equation (\ref{Back}), then
$u$ is given by (\ref{Rep}).
\et

\begin{proof}
First of all, let
\be
w_t(x):=\mE \Lambda^u_{u_0}(t,x)-\mE \Lambda^u_{f}(t,x).\label{De1}
\ee
By Theorem \ref{Th11}, $w(t,x)=w_t(x)$ solves the following backward Kolmogorov's equation:
\be
\p_t w+(u\cdot\nabla)w+(\nabla^\mathrm{t}u) w+\nu \Delta w=f,\ \ w(0,x)=u_0(x).\label{De2}
\ee
In view of $u=\bP w$, we may write
$$
w=u+\nabla q.
$$
Substituting it into  equation (\ref{De2}), one finds that
$$
\p_t u+(u\cdot\nabla)u+\nu \Delta u+\nabla p=f,\ \ u(0,x)=u_0(x),
$$
where
$$
p=\p_t q+(u\cdot\nabla)q+\nu\Delta q+\frac{1}{2}|u|^2.
$$

Conversely, let $(u,p)$ solve (\ref{Back}).
As above, if we define $w$ by (\ref{De1}), then $w$ satisfies equation (\ref{De2}).
We need to show that $u=\bP w$, or equivalently, for some scalar valued function $q$
$$
v:=w-u=\nabla q.
$$
By (\ref{De2}) and (\ref{Back}), $v$ solves the following equation:
\be
\p_t v+(u\cdot\nabla)v+(\nabla^\tr u)v+\nu \Delta v=\nabla p-\frac{1}{2}\nabla|u|^2,\ \ v(0,x)=0.\label{Es0}
\ee
Let
$$
q(t,x):=\mE\left(\int^0_t\Big[\frac{1}{2}|u(r,X_{t,r}(x))|^2-p(r,X_{t,r}(x))\Big]\dif r\right).
$$
Then $q$ solves the following equation (cf. \cite{Fr}):
$$
\p_t q+(u\cdot\nabla)q+\nu \Delta q=p-\frac{1}{2}|u|^2,\ \ q(0,x)=0.
$$
Taking gradients for both sides of the above equation yields
$$
\p_t \nabla q+(u\cdot\nabla)\nabla q+(\nabla^\tr u)(\nabla q)
+\nu \Delta \nabla q=\nabla p-\frac{1}{2}\nabla |u|^2.
$$
By the uniqueness of solutions to linear equation (\ref{Es0}), $v=\nabla q$.
\end{proof}
\br
In the above proof, we have assumed that the solutions are regular enough so that all the calculations
are valid. The existence of regular solutions will be proven in the next section.
\er

\section{A Proof of Local Existence in Sobolev Spaces}

With a little abuse of notations, in this and next sections we use $p$ to
denote the integrability index
since the pressure will not appear below.
For $k\in\mN_0$ and $p>1$, let $W^{k,p}(\mR^d;\mR^d)$ be the usual
$\mR^d$-valued Sobolev space on $\mR^d$, i.e,
the completion of $C^\infty_0(\mR^d;\mR^d)$ with respect to the norm:
$$
\|u\|_{k,p}:=\|u\|_{p}+\sum_{j=1}^k\|\nabla^j u\|_{p},
$$
where $\|\cdot\|_p$ is the usual $L^p$-norm, and $\nabla^j$ is the $j$-order gradient operator.
Note that $W^{0,p}(\mR^d;\mR^d)=L^p(\mR^d;\mR^d)$ and the following Sobolev's embedding holds: for $p>d$ (cf. \cite{Fr0})
\be
W^{1,p}(\mR^d;\mR^d)\hookrightarrow L^\infty(\mR^d;\mR^d),\ i.e. \
\|\cdot\|_\infty\leq c\|\cdot\|_{1,p},\label{Sob}
\ee
where $c=c(p,d)$. Below, we shall use $c$ to denote a constant  which may change in different occasions,
and whose dependence on  parameters can be traced carefully from the context.
Let $W^{k,p}_{loc}(\mR^d;\mR^d)$ be the local Sobolev space on $\mR^d$.
We introduce the following Banach space of transformations of $\mR^d$:
$$
\mX^{k+2,p}:=\Big\{X\in W^{k+2,p}_{loc}(\mR^d;\mR^d): |X(0)|
+\|\nabla X\|_\infty+\|\nabla^2 X\|_{k,p}<+\infty\Big\}.
$$

\bd
The Weber operator $\bW: L^p(\mR^d;\mR^d)\times\mX^{2,p}\to L^p(\mR^d;\mR^d)$ is defined by
$$
\bW(v,\ell):=\bP[(\nabla^\tr\ell)v],
$$
where $\bP$ is the Leray-Hodge projection onto divergence free vector fields.
\ed
\br
$\bP=I-\nabla(-\Delta)^{-1}\mathrm{div}$ is a singular integral operator(SIO) which is bounded
in $L^p$-space for $p\in(1,\infty)$(cf. \cite{St}).
\er

We now prepare several lemmas for later use.
\bl
(i) For any $k\in\mN_0$ and $p>d$, there exists a constant $c=c(k,p,d)>0$ such that
for all $v\in W^{k+2,p}(\mR^d;\mR^d)$ and $\ell\in \mX^{k+2,p}$,
\be
\|\nabla\bW(v,\ell)\|_{k+1,p}\leq c(\|\nabla\ell\|_{\infty}+\|\nabla^2\ell\|_{k,p})\cdot
\|\nabla v\|_{k+1,p}
\label{PL1}
\ee

(ii) For  $p>d$, there exists a constant $c=c(p,d)>0$ such that
for all $v_1,v_2\in W^{2,p}(\mR^d;\mR^d)$ and $\ell_1,\ell_2\in \mX^{2,p}$ with
$\ell_1-\ell_2\in L^p(\mR^d;\mR^d)$
\be
\|\bW(v_1,\ell_1)-\bW(v_2,\ell_2)\|_{p}\leq c
\Big(\|v_1\|_{2,p}\|\ell_1-\ell_2\|_{p}
+\|\nabla\ell_2\|_{\infty}\|v_1-v_2\|_{p}\Big).\label{PL2}
\ee
\el
\begin{proof}
(i)  Noting that
\be
\bP((\nabla^\tr \ell)v)+\bP((\nabla^\tr v)\ell)=\bP(\nabla(\ell\cdot v))=0,\label{Es5}
\ee
we have
$$
\p_i\bW(v,\ell)=\bP[(\nabla^\tr\p_i\ell)v+(\nabla^\tr\ell)\p_i v]
=\bP[-(\nabla^\tr v)\p_i\ell+(\nabla^\tr\ell)\p_i v]
$$
and
\ce
\p_j\p_i\bW(v,\ell)=\bP[-(\nabla^\tr v)\p_j\p_i\ell-(\nabla^\tr\p_jv)\p_i\ell
+(\nabla^\tr\ell)\p_j\p_i v+(\nabla^\tr\p_j\ell)\p_i v].
\de
Hence, by (\ref{Sob}) we have
\ce
\|\p_i\bW(v,\ell)\|_p&\leq& c(\|(\nabla^\tr v)\p_i\ell\|_p+\|(\nabla^\tr\ell)\p_i v\|_p)\\
&\leq& c(\|\nabla v\|_p\cdot\|\p_i\ell\|_\infty+\|\nabla\ell\|_\infty\|\p_i v\|_p)\\
&\leq&c\|\nabla v\|_p\cdot\|\nabla\ell\|_\infty
\de
and
\ce
\|\p_j\p_i\bW(v,\ell)\|_{p}&\leq& c(\|\nabla^2\ell\|_{p}
\cdot\|\nabla v\|_{\infty}+\|\nabla\ell\|_{\infty}\cdot\|\nabla^2v\|_{p})\\
&\leq& c(\|\nabla^2\ell\|_{p}
\cdot\|\nabla v\|_{1,p}+\|\nabla\ell\|_{\infty}\cdot\|\nabla^2v\|_{p})\\
&\leq&c(\|\nabla\ell\|_{\infty}+\|\nabla^2\ell\|_{p})\cdot\|\nabla v\|_{1,p},
\de
which produces
\ce
\|\nabla\bW(v,\ell)\|_{1,p}\leq c(\|\nabla\ell\|_{\infty}+\|\nabla^2\ell\|_{p})\cdot\|\nabla v\|_{1,p}.
\de
The higher derivatives can be estimated similarly.

(ii) By (\ref{Es5}), we have
\ce
\bW(v_1,\ell_1)-\bW(v_2,\ell_2)&=&\bP((\nabla^\tr(\ell_1-\ell_2))v_1)+\bP((\nabla^\tr\ell_2)(v_1-v_2))\\
&=&-\bP((\nabla^\tr v_1)(\ell_1-\ell_2))+\bP((\nabla^\tr\ell_2)(v_1-v_2)).
\de
So,
\ce
\|\bW(v_1,\ell_1)-\bW(v_2,\ell_2)\|_p&\leq&c\|(\nabla^\tr v_1)(\ell_1-\ell_2)\|_p
+c\|(\nabla^\tr\ell_2)(v_1-v_2)\|_p\\
&\leq&c\|\nabla v_1\|_\infty\|\ell_1-\ell_2\|_p
+c\|\nabla\ell_2\|_\infty\|v_1-v_2\|_p,
\de
which yields (\ref{PL2}) by (\ref{Sob}).
\end{proof}
\bl\label{Le1}
(i) For $k\in\mN_0$ and $p>d$, there exist  constants $c=c(k,p,d)>0$ and
$\a_k\in\mN_0$ such that for all $u\in W^{k+2,p}(\mR^d;\mR^d)$
and all  $X\in \mX^{k+2,p}$ preserving the volume,
\be
\|\nabla(u\circ X)\|_{k+1,p}\leq c\|\nabla u\|_{k+1,p}(1+\|\nabla X\|^{k+2}_{\infty}
+\|\nabla^2 X\|_{k,p}^{\a_k}).\label{PL3}
\ee
(ii) For $p>d$, there exists a constant $c=c(p,d)>0$ such that for all
$u\in W^{2,p}(\mR^d;\mR^d)$ and
$X,\tilde X\in\mX^{2,p}$ with $X-\tilde X\in L^p(\mR^d;\mR^d)$,
\be
\|u\circ X-u\circ \tilde X\|_{p}\leq
c\|\nabla u\|_{1,p} \cdot\|X-\tilde X\|_{p}.\label{PL33}
\ee
\el
\begin{proof}
(i) Since $X$ preserves the volume, we have
\be
\|u\circ X\|_{p}=\|u\|_{p}.\label{Es6}
\ee
Observe that for $m\geq 2$
\be
\nabla^k(u\circ X)=(\nabla^m u)\circ X\cdot(\nabla X)^m+\cdots
+(\nabla u)\circ X\cdot\nabla^{m} X.\label{Es8}
\ee
(\ref{PL3}) follows by (\ref{Sob}) and (\ref{Es6}).

(ii) It follows from
\be
u\circ X- u\circ\tilde X=\int^1_0(\nabla u)\circ (s X
+(1-s)\tilde X))\cdot(X-\tilde X)\dif s\label{PL22}
\ee
and (\ref{Sob}).
\end{proof}
\bl\label{Le2}
For $k\in\mN_0$, $U>0$ and $T:=\frac{1}{U}$, there exist
constants $c_1=c_1(p,d)>0$ and $c_2=c_2(k,p,d)>0$
such that for any divergence free vector field
$u\in C([-T,0];W^{k+2,p})$ satisfying
$\sup_{t\in[-T,0]}\|\nabla u_t\|_{k+1,p}\leq U$,
the solution $X_{t,s}$ to (\ref{Se1}) belongs to $\mX^{k+2,p}$ a.s., and
for all $t\in[-T,0]$
\be
\|\nabla X_{t,0}\|_{\infty}\leq  c_1,\ \
\|\nabla^2 X_{t,0}\|_{k,p}\leq  c_2.\label{PL5}
\ee
\el
\begin{proof}
Noting that
\be
\nabla X_{t,s}=\mI+\int^s_t(\nabla u_r)\circ X_{t,r}\cdot \nabla X_{t,r}\dif r,\label{Es77}
\ee
we have
$$
\|\nabla X_{t,s}\|_{\infty}\leq 1+\int^s_t\|\nabla X_{t,r}\|_{\infty}\cdot
\|\nabla u_r\|_{\infty}\dif r.
$$
By Gronwall's inequality,  we obtain  by (\ref{Sob})
\be
\sup_{-T\leq t\leq s\leq 0}\|\nabla X_{t,s}\|_{\infty}\leq
\exp\left[\int^0_{-T}\|\nabla u_r\|_{\infty}\dif r\right]\leq e^{cUT}=e^c=: c_1.\label{Es9}
\ee

On the other hand, from (\ref{Es77}), an elementary calculation shows that
$$
\mathrm{det}(\nabla X_{t,s})=\exp\left[\int^s_t(\nabla
\cdot u_r)\circ X_{t,r}\dif r\right]=1.
$$
So, $x\mapsto  X_{t,s}(x)$ preserves the volume.
Now,
$$
\nabla^2 X_{t,s}=\int^s_t[(\nabla^2 u_r)\circ X_{t,r}\cdot (\nabla X_{t,r})^2
+(\nabla u_r)\circ X_{t,r}\cdot \nabla^2 X_{t,r}]
\dif r.
$$
Hence, by (\ref{Es9}) and (\ref{Sob})
\ce
\|\nabla^2 X_{t,s}\|_{p}&\leq&c\int^s_t\Big(\|\nabla^2u_r\|_{p}
\|\nabla X_{t,r}\|_{\infty}^2+\|\nabla u_r\|_{\infty}\|\nabla^2X_{t,r}\|_{p}\Big)\dif r\\
&\leq&cUT+cU\int^s_t\|\nabla^2X_{t,r}\|_{p}\dif r.
\de
By Gronwall's inequality again  we get
$$
\|\nabla^2 X_{t,s}\|_{p}\leq cUTe^{cUT}=ce^{c}=:c_2.
$$
Higher derivatives can be estimated similarly step by step.
\end{proof}
\bl\label{Le3}
For $p>d$ and $T>0$, let $u,\tilde u\in C([-T,0]; W^{2,p}(\mR^d;\mR^d))$, and
$X,\tilde X$ solve SDE (\ref{SDE}) with drifts
$u$ and $\tilde u$ respectively.
Then for some $c=c(p,d)>0$ and any $t\in[-T,0]$,
\be
\|X_{t,0}-\tilde X_{t,0}\|_{p}\leq \exp\Big[c\sup_{t\in[-T,0]}\|\nabla u_t\|_{1,p}\Big]
\cdot\int^0_t\|u_r-\tilde u_r\|_{p}\dif r.\label{PL6}
\ee
\el
\begin{proof}
We have
\ce
X_{t,s}(x)-\tilde X_{t,s}(x)&=&\int^s_t(u_r(X_{t,r}(x))-\tilde u_r(\tilde X_{t,r}(x)))\dif r\\
&=&\int^s_t(u_r(X_{t,r}(x))-u_r(\tilde X_{t,r}(x)))\dif r\\
&&+\int^s_t(u_r(\tilde X_{t,r}(x))-\tilde u_r(\tilde X_{t,r}(x)))\dif r.
\de
For $R>0$, let $B_R:=\{x\in\mR^d: |x|\leq R\}$. By virtue of $x\mapsto\tilde X_{t,r}(x)$
preserving the volume and formula (\ref{PL22}), we have
\be
\|X_{t,s}-\tilde X_{t,s}\|_{L^p(B_R)}&\leq&
\int^s_t\|u_r\circ X_{t,r}- u_r\circ\tilde X_{t,r}\|_{L^p(B_R)}\dif r
+\int^s_t\|u_r-\tilde u_r\|_{p}\dif r\no\\
&\leq&\sup_{r\in[-T,0]}\|\nabla u_r\|_\infty\int^s_t\|X_{t,r}-\tilde X_{t,r}\|_{L^p(B_R)}\dif r\no\\
&&+\int^0_t\|u_r-\tilde u_r\|_{p}\dif r,\label{Es4}
\ee
By Gronwall's inequality  and (\ref{Sob}), we get
$$
\|X_{t,s}-\tilde X_{t,s}\|_{L^p(B_R)}\leq \exp\Big[c\sup_{t\in[-T,0]}\|\nabla u_t\|_{1,p}\Big]
\cdot\int^0_t\|u_r-\tilde u_r\|_{p}\dif r.
$$
Letting $R$ go to infinity gives (\ref{PL6}).
\end{proof}

We are now in a position to prove the following local existence result.
\bt\label{Main}
For $\nu\geq 0$, $k\in\mN_0$ and $p>d$,
there exists a constant $c_0=c_0(k,p,d)>0$ independent of $\nu$ such that for
any $u_0\in W^{k+2,p}(\mR^d;\mR^d)$ divergence free and $T:=(c_0\|\nabla u_0\|_{k+1,p})^{-1}$,
there is a unique pair $(u,X)$ with
$u\in C([-T,0];W^{k+2,p})$ satisfying
\be\label{SN}
\left\{
\begin{aligned}
&X_{t,s}(x)=x+\int^s_tu_r(X_{t,r}(x))\dif r+\sqrt{2\nu}\int^s_t\sigma_r\dif B_r,\
t\leq s\leq 0,\\
&u_t=\bP\mE [(\nabla^\tr X_{t,0})(u_0\circ X_{t,0})],\ \ t\leq0.
\end{aligned}
\right.
\ee
Moreover, for any $t\in[-T,0]$
\be
\|\nabla u_t\|_{k+1,p}\leq c_0\|\nabla u_0\|_{k+1,p}.\label{PL7}
\ee
\et
\begin{proof}
Set $u^1_r(x):=u_0(x)$. Consider the following Picard's iteration sequence
\be\label{App}
\left\{
\begin{aligned}
X^n_{t,s}(x)&=x+\int^s_tu^n_r(X^n_{t,r}(x))\dif r+\sqrt{2\nu}\int^s_t\sigma(r)\dif B_r,\ t\leq s\leq 0,\\
u^{n+1}_t&=\bP\mE [(\nabla^\tr X^n_{t,0})(u_0\circ X^n_{t,0})],\ \ t\leq0.
\end{aligned}
\right.
\ee
Noting that
$$
\bP\mE [(\nabla^\tr X^n_{t,0})(u_0\circ X^n_{t,0})]=\mE\bW(u_0\circ X^n_{t,0},X^n_{t,0}),
$$
we have by (\ref{PL1}) and (\ref{PL3})
\ce
\|\nabla u^{n+1}_t\|_{k+1,p}&\leq& \mE\|\nabla\bW(u_0\circ X^n_{t,0},X^n_{t,0})\|_{k+1,p}\\
&\leq&c\mE\Big[(\|\nabla X^n_{t,0}\|_{\infty}+\|\nabla^2 X^n_{t,0}\|_{k,p})
\cdot\|\nabla(u_0\circ X^n_{t,0})\|_{k+1,p}\Big]\\
&\leq&c_3\mE\Big[(1+\|\nabla X^n_{t,0}\|^{k+3}_{\infty}+\|\nabla^2 X^n_{t,0}\|^{\beta_k}_{k,p})
\cdot\|\nabla u_0\|_{k+1,p}\Big],
\de
where $\beta_k\in\mN$ only depends on $k$ and $c_3=c_3(k,p,d)\geq 1$.

Set
$$
c_0:=c_3(1+c_1^{k+3}+c_2^{\beta_k})\geq 1,
$$
where $c_1$ and $c_2$ are from Lemma \ref{Le2}.
Choosing $U:=c_0\|\nabla u_0\|_{k+1,p}$ and $T:=1/U$ in Lemma \ref{Le2},
we have by induction and Lemma \ref{Le2}
\be
\sup_{t\in[-T,0]}\|\nabla u^{n}_t\|_{k+1,p}\leq U,\ \ \forall n\in\mN.\label{Es10}
\ee
On the other hand, we also have by (\ref{PL5})
\ce
\|u^{n+1}_t\|_{p}&\leq&c\mE\|(\nabla^\tr X^n_{t,0})(u_0\circ X^n_{t,0})\|_p\\
&\leq&c\mE\Big[\|\nabla^\tr X^n_{t,0}\|_\infty\|u_0\circ X^n_{t,0}\|_p\Big]\\
&\leq&c\|u_0\|_p,
\de
which together with (\ref{Es10}) gives the following uniform estimate:
\be
\sup_{n\in\mN}\sup_{t\in[-T,0]}\|u^{n}_t\|_{k+2,p}<+\infty.\label{Es11}
\ee

Now by (\ref{PL2}) and (\ref{PL33}) (\ref{PL6}), we have
\ce
\|u^{n+1}_t-u^{m+1}_t\|_{p}&\leq& c\mE\Big[\|u_0\circ X^n_{t,0}\|_{2,p}
\cdot\| X^n_{t,0}-X^m_{t,0}\|_{p}\\
&&+\|\nabla X^m_{t,0}\|_\infty\cdot\|u_0\circ X^n_{t,0}-u_0\circ X^m_{t,0}\|_{p}\Big]\\
&\leq& c\int^0_t\|u^n_r- u^m_r\|_{p}\dif r,
\de
where $c=c(p,d,U)$ is independent of $n,m$. From this we derive that
$$
\limsup_{n,m\to\infty}\sup_{t\in[-T,0]}\|u^n_t- u^m_t\|_{p}=0.
$$
By (\ref{Es11}) and interpolation inequality, we further have
$$
\limsup_{n,m\to\infty}\sup_{t\in[-T,0]}\|u^n_t- u^m_t\|_{k+1,p}=0.
$$
Therefore, there is a $u\in C([-T,0];W^{k+1,p}(\mR^d;\mR^d))$ such that
$$
\limsup_{n\to\infty}\sup_{t\in[-T,0]}\|u^n_t- u_t\|_{k+1,p}=0.
$$
Taking limits for (\ref{App}), one finds that $u$ is a solution of (\ref{SN}).
Estimate (\ref{PL7}) follows from (\ref{Es10}).
\end{proof}
\br
The constant $c_0$ in (\ref{PL7}) is usually strictly greater than $1$. If $c_0$ equals $1$, then we
can invoke the standard bootstrap method to obtain the global existence. This will be studied in
the next section when the periodic boundary is considered and the
viscosity is large enough.
\er

Since the existence time interval in Theorem \ref{Main} is independent of
the viscosity $\nu$, we also obtain the local existence of solutions to Euler equation (\ref{Euler}).
Moreover, as $\nu\to 0$, the solution of Navier-Stokes equation
converges to the solution of Euler equation as given below.
\bp\label{Pr}
Keep the same assumptions as in Theorem \ref{Main}. For $\nu\geq 0$ and $u_0\in W^{k+2,p}\cap L^2$,
let $(u^\nu, X^\nu)$ be the solution of (\ref{SN}) corresponding to viscosity $\nu$ and initial value
$u_0$.
Then for any $j=0,\cdots,k+1$, there exists $c=c(k,j,p,d,\|u_0\|_{k+2,p},\|u_0\|_2)>0$
such that for all $\nu\geq 0$ and $t\in[-T,0]$
$$
\|u^\nu_t-u^0_t\|_{C^j_b}\leq c(\nu |t|)^{(\frac{k+2-j}{d}-\frac{1}{p})
/(\frac{1}{2}+\frac{k+2}{d}-\frac{1}{p})}.
$$
\ep
\begin{proof}
Note that $u_0\in W^{k+2,p}\cap L^2$ guarantees $u^\nu_t\in W^{k+2,p}\cap L^2$. By
$$
\p_t(u^\nu_t-u^0_t)+\nu u^\nu_t+\bP[(u^\nu_t\cdot\nabla)u^\nu_t-(u^0_t\cdot\nabla)u^0_t]=0,
$$
we have
\ce
-\p_t\|u^\nu_t-u^0_t\|_2^2&=&\nu\<\Delta u^\nu_t,u^\nu_t-u^0_t\>_2
+\<(u^\nu_t\cdot\nabla)u^\nu_t-(u^0_t\cdot\nabla)u^0_t,u^\nu_t-u^0_t\>_2\\
&=&\nu\<\Delta u^\nu_t,u^\nu_t-u^0_t\>_2
+\<((u^\nu_t-u^0_t)\cdot\nabla)u^\nu_t,u^\nu_t-u^0_t\>_2\\
&\leq&\nu\|\Delta u^\nu_t\|_2\cdot\|u^\nu_t-u^0_t\|_2
+\|\nabla u^\nu_t\|_\infty\|u^\nu_t-u^0_t\|^2_2,
\de
i.e.,
$$
-\p_t\|u^\nu_t-u^0_t\|_2\leq\nu\|\Delta u^\nu_t\|_2
+\|\nabla u^\nu_t\|_\infty\|u^\nu_t-u^0_t\|_2.
$$
By Gronwall's inequality and (\ref{PL7}) we obtain
\ce
\|u^\nu_t-u^0_t\|_2\leq\nu\int^0_t\|\Delta u^\nu_s\|_2\dif s\cdot
\exp\left[\int^0_t\|\nabla u^\nu_s\|_\infty\dif s\right]\leq c\nu |t|.
\de
The desired estimate now follows by the Sobolev embedding(cf. \cite{Fr0}): for $u\in W^{k+2,p}\cap L^2$
$$
\|\nabla^j u\|_\infty\leq c_{k,j,p,d}\|u\|^\a_{k+2,p}\|u\|_2^{1-\a},
$$
where $\a=(\frac{j}{d}+\frac{1}{2})/(\frac{1}{2}+\frac{k+2}{d}-\frac{1}{p})$.
\end{proof}
\br
We cannot prove a convergence rate $O(\sqrt{\nu t})$ as in \cite{Iy} starting from
(\ref{SN}) because $x\mapsto (X^\nu_{t,0}(x)-X^0_{t,0}(x))$ does not belong to any $L^p$-spaces.
\er

\section{Existence of Global Solutions}

\subsection{Global Existence in Two Dimensions}
First of all, we recall the following Beale-Kato-Majda's estimate about
SIOs, which can be proved as in \cite[p.117, Proposition 3.8]{Be-Ma},
we omit the details.
\bl
For $p>d$, let $u\in W^{2,p}(\mR^d;\mR^d)$ be a divergence free vector field and
$\omega:=\mathrm{curl}u$. Then, for some $c=c(p,d)$
\be
\|\nabla u\|_\infty\leq c(1+\log^+\|\omega\|_{1,p})(1+\|\omega\|_\infty),\label{DW5}
\ee
where $\log^+x:=\max\{\log x, 0\}$ for $x>0$.
\el
In two dimensional case,  taking the curl for the second equation in (\ref{SN}), one finds that
\be
\omega_t:=\mathrm{curl} u_t:=\p_1 u^2_t-\p_2 u^1_t
=\mE [\omega_0\circ X_{t,0}].\label{DW4}
\ee
From this, we clearly have
\be
\|\omega_t\|_p\leq\|\omega_0\|_p,\ \ 1\leq p\leq \infty.\label{DW1}
\ee
Basing (\ref{DW5}) and representation (\ref{DW4}), we may prove the following global existence
for 2D Navier-Stokes and Euler equations.
\bt\label{Th1}
In two dimensions, for $\nu\geq 0$, $k\in \mN_0$, $p>2$ and
$u_0\in W^{k+2,p}(\mR^2;\mR^2)$ divergence free,
there exists a unique global solution $(u,X)$ to equation (\ref{SN}).
\et
\begin{proof}
We only need to prove the following a priori estimate: for all $t\in\mR_-$
$$
\|u_t\|_{k+2,p}\leq c(\|u_0\|_{k+2,p},k,p,t)<+\infty,
$$
where $c(\|u_0\|_{k+2,p},k,p,t)$ continuously depends on its parameters.

Following the proof of Lemma \ref{Le2}, we have
\be
\|\nabla X_{t,0}\|_\infty\leq\exp\left[\int^0_{t}\|\nabla u_r\|_{\infty}\dif r\right].\label{DW2}
\ee
Noting that
$$
\nabla\omega_t=\mE(\nabla\omega_0\circ X_{t,0}\cdot\nabla X_{t,0}),
$$
we have
$$
\|\nabla\omega_t\|_p\leq\|\nabla\omega_0\|_p\cdot\mE\|\nabla X_{t,0}\|_\infty
$$
and by (\ref{DW1}) and (\ref{DW2})
\be
\|\omega_t\|_{1,p}\leq\|\omega_0\|_{1,p}\cdot\left(1+
\exp\left[\int^0_{t}\|\nabla u_r\|_{\infty}\dif r\right]\right).\label{DW3}
\ee
Hence, by (\ref{DW5}) (\ref{DW1}) and (\ref{DW3})
\ce
\|\nabla u_t\|_\infty&\leq& c(1+\log^+\|\omega_t\|_{1,p})(1+\|\omega_t\|_\infty)\\
&\leq&c +c\int^0_{t}\|\nabla u_r\|_{\infty}\dif r,
\de
where $c=c(\|\omega_0\|_{1,p},p)$.
By Gronwall's inequality we obtain
$$
\|\nabla u_t\|_\infty\leq ce^{c|t|}.
$$
Substituting this into (\ref{DW2}) and (\ref{DW3}) gives
$$
\|\nabla X_{t,0}\|_\infty\leq e^{c|t| e^{c|t|}},
$$
and by Calderon-Zygmund's inequality about SIOs (cf. \cite{St})
$$
\|\nabla u_t\|_{1,p}\leq\|\omega_t\|_{1,p}\leq\|\omega_0\|_{1,p}\cdot\left(1+e^{c|t| e^{c|t|}}\right).
$$
Moreover,
$$
\|u_t\|_p\leq c\mE(\|\nabla X_{t,0}\|_\infty\cdot\|u_0\circ X_{t,0}\|_p)
\leq c\|u_0\|_p\cdot e^{c|t| e^{c|t|}}.
$$
Thus,
$$
\|u_t\|_{2,p}\leq c(\|u_0\|_{2,p},p,t)<+\infty.
$$
Starting from (\ref{DW4}) and as in Lemma \ref{Le2}, higher derivatives can be estimated
similarly.
\end{proof}
\subsection{Global Existence for Large Viscosity}
In this section, we study the existence of global solutions for large viscosity and work
on the $d$-dimensional torus $\mT^d=\mR^d/\mZ^d$.
Let $W^{k,p}(\mT^d,\mR^d)$ be the $\mR^d$-valued Sobolev spaces
on $\mT^d$ with vanishing mean.
Instead of (\ref{Se1}), we consider
\be
X_{t,s}(x)=x+\int^s_tu_r(X_{t,r}(x))\dif r+\sqrt{2\nu}(B_s-B_t),\label{Lp2}
\ee
where $B$ is the standard Wiener process on $\Omega:=C(\mR_-;\mR^d)$, i.e., for $\omega\in\Omega$,
$B_\cdot(\omega)=\omega(\cdot)$.

We recall the following Bismut's formula (cf. \cite{Bi, Cr-Zh}).
For the reader's convenience, a short proof
is provided here.
\bt
For any $t<0$ and $f\in C_b^1(\mT^d;\mR)$, it holds that
\be
(\nabla\mE f(X_{t,0}))(x)=\frac{1}{t\sqrt{2\nu}}\mE\left[f(X_{t,0}(x))
\int^0_t\Big(s(\nabla u_s)\circ X_{t,s}(x)-\mI\Big)\dif B_s\right].\label{For}
\ee
In particular, for any $p>d$ and some $c=c(p,d)$
\be
\|\nabla\mE f(X_{t,0})\|_{p}\leq\frac{c}{\sqrt{\nu |t|}}\|f\|_p
\left[|t|\cdot\sup_{s\in[t,0]}\|\nabla u_s\|_{1,p}+1\right].\label{Lp3}
\ee
\et
\begin{proof}
Fix $t<0$ and $y\in \mR^d$ below and define
$$
h(s):=\frac{1}{t\sqrt{2\nu}}\left[(t-s)y
+\int^s_t[(\nabla u_r)\circ X_{t,r}(x)]\cdot (ry)\dif r\right],\ s\in[t,0].
$$
Consider the Malliavin derivative of $X_{t,s}$ with respect to the sample path along
the direction $h$, i.e.,
$$
D_h X_{t,s}(x,\omega)=\lim_{\eps\to 0}\frac{X_{t,s}(x,\eps h+\omega)-X_{t,s}(x,\omega)}{\eps},
\ \omega\in\Omega.
$$
From (\ref{Lp2}) one sees that
\ce
D_h X_{t,s}(x)&=&\int^s_t [(\nabla u_r)\circ X_{t,r}(x)]\cdot D_h X_{t,r}(x)\dif r+\sqrt{2\nu}h(s)\\
&=&\frac{(t-s)y}{t}+\int^s_t [(\nabla u_r)\circ X_{t,r}(x)]\cdot
\Big[D_h X_{t,r}(x)+\frac{ry}{t}\Big]\dif r.
\de
On the other hand,  we have
$$
\nabla X_{t,s}(x)\cdot y
=y+\int^s_t [(\nabla u_r)\circ X_{t,r}(x)]\cdot \nabla X_{t,r}(x)\cdot y\dif r.
$$
By the uniqueness of solutions, we get
$$
\nabla X_{t,s}(x)\cdot y=D_h X_{t,s}(x)+\frac{sy}{t}.
$$
In particular,
$$
\nabla X_{t,0}(x)\cdot y=D_h X_{t,0}(x).
$$

Now
\ce
\nabla\mE f(X_{t,0})\cdot y&=&\mE\big[[(\nabla f)\circ X_{t,0}]\cdot\nabla X_{t,0}\cdot y\big]\\
&=&\mE\big[(\nabla f)\circ X_{t,0}\cdot D_h X_{t,0}\big]\\
&=&\mE\big[D_h(f\circ X_{t,0})\big]\\
&=&\mE\left[(f\circ X_{t,0})\int^0_t\dot h(s)\dif B_s\right],
\de
where the last step is due to the integration by parts formula in the Malliavin calculus
(cf. \cite{Ma}). Formula (\ref{For}) now follows.

For estimation (\ref{Lp3}), by H\"older's inequality and It\^o's isometry,
the square of the right hand side of (\ref{For}) is controlled by
\ce
&&\frac{1}{2\nu t^2}\mE|f(X_{t,0}(x))|^2\mE\left[\int^0_t
|s(\nabla u_s)\circ X_{t,s}(x)-\mI|^2\dif s\right]\\
&&\qquad\leq\frac{c}{\nu t^2}\mE|f(X_{t,0}(x))|^2
\left[|t|^3\sup_{s\in[t,0]}\|\nabla u_s\|^2_\infty+|t|\right].
\de
Hence, by (\ref{Sob})
\ce
\|\nabla\mE f(X_{t,0})\|_{p}&\leq&\frac{c}{\sqrt{\nu |t|}}\|f\|_p
\left[|t|\cdot\sup_{s\in[t,0]}\|\nabla u_s\|_\infty+1\right]\\
&\leq&\frac{c}{\sqrt{\nu |t|}}\|f\|_p
\left[|t|\cdot\sup_{s\in[t,0]}\|\nabla u_s\|_{1,p}+1\right].
\de
The proof is complete.
\end{proof}

We now prove the following global existence result (see also \cite{Ko-Ta, Iy0}).
\bt
Let $k\in\mN_0$ and $p>d$, $u_0\in W^{k,p}(\mT^d;\mR^d)$ be divergence free and mean zero.
Let $(u,X)$ be the local solution of (\ref{SN}) in Theorem \ref{Main}.
Then, there exist $T_0=T_0(k,p,d,\|\nabla u_0\|_{k+1,p})<0$ and $\d=\d(k,p)>0$
such that if $\nu\geq\d\|\nabla u_0\|_{k+1,p}$, then
$$
\|\nabla u_{T_0}\|_{k+1,p}\leq \|\nabla u_0\|_{k+1,p},
$$
and there is a global solution to equation (\ref{SN}).
\et
\begin{proof}
Let $(u,X)$ be the local solution  of (\ref{SN}) on $[-T,0]$ in Theorem \ref{Main}, where
$T=(c_0\|\nabla u_0\|_{k+1,p})^{-1}$.
Recalling the estimations in Lemma \ref{Le2} and Theorem \ref{Main},  we have
for all  $t\in[-T,0]$
\be
\|\nabla X_{t,0}\|_\infty\leq c_1,\ \ \|\nabla^2 X_{t,0}\|_{k,p}\leq c_2\label{PL99}
\ee
and
\be
\|\nabla u_t\|_{k+1,p}\leq c_0\|\nabla u_0\|_{k+1,p}.\label{PL9}
\ee

Write
\be
u_t=\bP\mE[(\nabla^\tr X_{t,0}-\mI)(u_0\circ X_{t,0})]+\bP\mE(u_0\circ X_{t,0}).\label{PL8}
\ee
We separately deal with the first term and the second term.

For the first term in (\ref{PL8}), using (\ref{PL99}) (\ref{PL9}) and as in Lemma \ref{Le2},
one may prove that for some $c=c(k,p,d)$ and
all $t\in[-T,0]$
$$
\|\nabla^\tr X_{t,0}-\mI\|_{k+2,p}\leq c\|\nabla u_0\|_{k+1,p}\cdot |t|.
$$
Using this estimate as well as (\ref{PL3}) (\ref{PL99}) and (\ref{PL9}), one finds that
\be
\|\nabla\bP\mE[(\nabla^\tr X_{t,0}-\mI)(u_0\circ X_{t,0})]\|_{k+1,p}
\leq c\|\nabla u_0\|_{k+1,p}^2\cdot |t|.\label{Op1}
\ee

For the second term in (\ref{PL8}), by  (\ref{Lp3}), (\ref{PL9})
and Poincare's inequality, we have
\be
\|\nabla\bP\mE(u_0\circ X_{t,0})\|_p\leq\frac{c}{\sqrt{\nu|t|}}
\|u_0\|_p\leq\frac{c}{\sqrt{\nu|t|}}
\|\nabla u_0\|_p.\label{Lp4}
\ee
Note that
$$
\nabla^2\mE(u_0\circ X_{t,0})=
\nabla\mE((\nabla u_0)\circ X_{t,0})+\nabla\mE[((\nabla u_0)\circ X_{t,0})(\nabla X_{t,0}-\mI)].
$$
As above, we have
$$
\|\nabla\mE((\nabla u_0)\circ X_{t,0})\|_p\leq\frac{c}{\sqrt{\nu |t|}}\|\nabla u_0\|_p
$$
and
$$
\|\nabla\mE[((\nabla u_0)\circ X_{t,0})(\nabla X_{t,0}-\mI)]\|_{k,p}
\leq  c\|\nabla u_0\|^2_{k+1,p}\cdot |t|.
$$
So,
\be
\|\nabla^2\bP\mE(u_0\circ X_{t,0})\|_p&\leq&c\|\nabla^2\mE(u_0\circ X_{t,0})\|_p\no\\
&\leq&\frac{c}{\sqrt{\nu |t|}}\|\nabla u_0\|_p
+c\|\nabla u_0\|^2_{k+1,p}\cdot |t|.\label{Lp5}
\ee
Continuing the above calculations we get
\be
\|\nabla^{k+2}\bP\mE(u_0\circ X_{t,0})\|_p\leq
\frac{c}{\sqrt{\nu |t|}}\|\nabla^{k+1}u_0\|_{p}
+c\|\nabla u_0\|_{k+1,p}^2\cdot |t|.\label{Lp6}
\ee
Combining (\ref{Lp4}) (\ref{Lp5}) and (\ref{Lp6}), we find
\be
\|\nabla\bP\mE(u_0\circ X_{t,0})\|_{k+1,p}\leq\frac{c}{\sqrt{\nu |t|}}\|\nabla u_0\|_{k,p}
+c\|\nabla u_0\|_{k+1,p}^2\cdot |t|.\label{Lp7}
\ee

Summarizing (\ref{PL8}) (\ref{Op1}) and (\ref{Lp7}) yields
$$
\|\nabla u_t\|_{k+1,p}\leq\left[\frac{c_3}{\sqrt{\nu |t|}}
+c_4\|\nabla u_0\|_{k+1,p}\cdot |t|\right]\|\nabla u_0\|_{k+1,p},
$$
where $c_3=c_3(k,p,d)$ and $c_4=c_4(k,p,d)>c_0$.
Now, taking $T_0=-\frac{1}{2c_4\|\nabla u_0\|_{k+1,p}}$ and $\delta=8c_3^2c_4$, we have
for $\nu\geq\delta\|\nabla u_0\|_{k+1,p}$
$$
\|\nabla u_{T_0}\|_{k+1,p}\leq \|\nabla u_0\|_{k+1,p}.
$$
The proof is thus finished.
\end{proof}

\section{A Large Deviation Estimate for Stochastic Particle Paths}

Let $\mG^k$ denote the  $k$-order diffeomorphism group on $\mR^d$, which is endowed with
the locally uniform convergence topology together with its inverse for all derivatives up to $k$.
Then $\mG^k$ is a Polish space.
Let $\mG^k_0$ be the subspace of $\mG^k$ in which each transformation preserves the
Lebesgue measure, equivalently,
$$
\mG^k_0:=\{X\in \mG^k: \mathrm{det}(\nabla X)=1\}.
$$
Then $\mG^k_0$ is a closed subspace of $\mG^k$, therefore, a Polish space.

It is clear that
$$
t\mapsto X^\nu_t(\cdot)\in\mG^k_0
$$
is continuous by the theory of stochastic flow (cf. \cite{Ku}).
We now state a large deviation principle of Freidlin-Wentzell's type,  which follows
from the results in \cite{Ben, RZ} by using Proposition \ref{Pr}.
\bt\label{t7}
Keep all the things as in Proposition \ref{Pr}.
For any Borel set $E\subset C([-T,0];\mG^k_0)$, we have
\ce
-\inf_{Y\in E^o}I(Y)\leq\liminf_{\nu\rightarrow 0}\nu\log\mP(X^\nu\in E)
\leq\limsup_{\nu\rightarrow 0}\nu\log\mP(X^\nu\in E)\leq -\inf_{Y\in \bar E}I(Y),
\de
where $E^o$ and $\bar E$ denotes the interior and the closure respectively in
$C([-T,0];\mG^k_0)$, and $I(Y)$ is the rate function defined by
$$
I(Y):=\frac{1}{2}\inf_{\{h\in L^2(-T,0;l^2): S(h)=Y\}}\int^0_{-T}\|h_s\|^2_{l^2}\dif s, \quad
Y\in C([-T,0];\mG^k_0),
$$
where $S(h)=Y$ solves the following ODE:
$$
Y_s(x)=x+\int^s_{-T} u^0_r(Y_r(x))\dif r+\int^s_{-T}\<\sigma_r, h_r\>_{l^2}\dif r,\ \ s\in[-T,0].
$$
\et
\br
In two dimensions, the $T$ in the above theorem
can be arbitrarily large by Theorem \ref{Th1}.
\er

\vspace{5mm}

{\bf Acknowledgements:}

The author would like to thank Professor Benjamin Goldys for
providing him an excellent environment to work in the University of New South Wales.
His work is supported by ARC Discovery grant DP0663153 of Australia.

\end{document}